\documentclass[psamsfonts,reqno]{amsart}
\usepackage{amssymb,eucal}

\usepackage{amscd}

\usepackage{graphics}

\usepackage{mathrsfs}

\newcommand{\color}[6]{}

\theoremstyle{plain}

\newtheorem{lemma}{Lemma}[section]
\newtheorem{theorem}[lemma]{Theorem}
\newtheorem{proposition}[lemma]{Proposition}
\newtheorem{corollary}[lemma]{Corollary}

\newtheorem*{stat}{\name}
\newcommand{\name}{testing}

\theoremstyle{definition}
\newtheorem{definition}[lemma]{Definition}

\theoremstyle{remark}
\newtheorem{remark}[lemma]{Remark}

\newcommand{\qedc}{{\qed}~{\rm Claim~{\theclaim}.}}
\newcommand{\qedsc}{{\qed}~{\rm Claim.}}

\numberwithin{equation}{section}

\newcommand{\mirr}{meet-ir\-re\-duc\-i\-ble}

\newcommand{\set}[1]{\{#1\}}
\newcommand{\setm}[2]{\set{#1\mid#2}}
\newcommand{\Set}[1]{\left\{#1\right\}}
\newcommand{\Setm}[2]{\Set{#1\mid#2}}

\newcommand{\famm}[2]{(#1\mid#2)}

\DeclareMathOperator{\im}{im}

\newcommand{\Pow}{\mathfrak{P}}

\newcommand{\sP}{\mathscr{P}}

\newcommand{\mbZ}{\mathbb{Z}}
\newcommand{\Powf}[1]{[#1]^{<\omega}}
\newcommand{\dnw}{\mathbin{\downarrow}}
\newcommand{\upw}{\mathbin{\uparrow}}

\newcommand{\cS}{\mathcal{S}}

\newcommand{\cV}{\mathcal{V}}
\newcommand{\cW}{\mathcal{W}}

\newcommand{\cM}{\mathcal{M}}

\newcommand{\bu}{\boldsymbol{u}}
\newcommand{\bv}{\boldsymbol{v}}

\DeclareMathOperator{\dom}{dom}
\DeclareMathOperator{\Con}{Con}
\DeclareMathOperator{\Conc}{Con_c}

\DeclareMathOperator{\Id}{Id}
\DeclareMathOperator{\Dim}{Dim}
\newcommand{\GDim}{K_0^\ell}

\newcommand{\id}{\mathrm{id}}
\newcommand{\jz}{$(\vee,0)$}

\newcommand{\jzs}{\jz-semi\-lat\-tice}

\newcommand{\jzh}{\jz-ho\-mo\-mor\-phism}

\newcommand{\res}{\mathbin{\restriction}}

\DeclareMathOperator{\crita}{crit}
\newcommand{\crit}[2]{\crita({{#1};{#2}})}

\newcommand{\module}[1]{|{#1}|}
\DeclareMathOperator{\card}{card}
\DeclareMathOperator{\Cond}{Cond}

\DeclareMathOperator{\M}{M}

\DeclareMathOperator{\Sub}{Sub}
\newcommand{\Var}{{\bf{Var}}}

\DeclareMathOperator{\Ide}{Id_e}
\DeclareMathOperator{\Ids}{Id_s}
\DeclareMathOperator{\Idc}{Id_c}

\DeclareMathOperator{\lh}{lh}

\newcommand{\FF}{\mathbb{F}}

\newcommand{\tosurj}{\mathbin{\twoheadrightarrow}}
\newcommand{\toinj}{\mathbin{\hookrightarrow}}

\newcommand{\un}{\underline{n}}

\newcommand{\mtol}{\nabla}
\newcommand{\pgrtol}{\overline{\nabla}}

\newcommand{\rideal}{\mathbb{L}}

\usepackage[all]{xy}

\begin{document}

\title[Critical points]{Critical points between varieties generated by subspace lattices of vector spaces}

\author[P.~Gillibert]{Pierre Gillibert}
\address{LMNO, CNRS UMR 6139\\
D\'epartement de Math\'ematiques, BP 5186\\
Universit\'e de Caen, Campus 2\\
14032 Caen cedex\\
France}
\address{
Charles University in Prague, Faculty of Mathematics and Physics, Department of Algebra, Sokolovska 83, 186 00 Prague, Czech Republic.
}
\email{pierre.gillibert@math.unicaen.fr, pgillibert@yahoo.fr}
\urladdr{http://www.math.unicaen.fr/\~{}giliberp/}



\thanks{This paper is a part of the author's ``Doctorat de l'universit\'e de Caen'', prepared under the supervision of Friedrich Wehrung}

\date{\today}

\begin{abstract}
We denote by $\Conc A$ the semilattice of all compact congruences of an algebra~$A$. Given a variety~$\cV$ of algebras, we denote by $\Conc\cV$ the class of all semilattices isomorphic to $\Conc A$ for some $A\in\cV$. Given varieties~$\cV$ and~$\cW$ of algebras, the \emph{critical point} of~$\cV$ under~$\cW$ is defined as $\crit{\cV}{\cW}=\min\setm{\card D}{D\in\Conc\cV - \Conc\cW}$. Given a finitely generated variety~$\cV$ of modular lattices, we obtain an integer~$\ell$, depending on~$\cV$, such that $\crit{\cV}{\Var(\Sub F^n)}\ge\aleph_2$ for any $n\ge\ell$ and any field~$F$.

In a second part, using tools introduced in \cite{G:critpoint}, we prove that:
\[
\crit{\cM_n}{\Var(\Sub F^3)}=\aleph_2,
\]
for any finite field~$F$ and any ordinal~$n$ such that $2+\card F\le n\le\omega$. Similarly $\crit{\Var(\Sub F^3)}{\Var(\Sub K^3)}=\aleph_2$, for all finite fields~$F$ and~$K$ such that $\card F>\card K$.
\end{abstract}

\maketitle

\section{Introduction}
We denote by~$\Con A$ (resp., $\Conc A$) the lattice (resp., \jzs) of all congruences (resp., compact congruences) of an algebra~$A$. For a homomorphism~$f\colon A\to B$ of algebras, we denote by $\Con f$ the map from~$\Con A$ to~$\Con B$ defined by the rule
 \[
 (\Con f)(\alpha)=\text{congruence of }B\text{ generated by }\setm{(f(x),f(y))}{(x,y)\in \alpha},
 \]
for every $\alpha\in\Con A$, and we also denote by~$\Conc f$ the restriction of~$\Con f$ from~$\Conc A$ to~$\Conc B$.

A \emph{congruence-lifting} of a \jzs~$S$ is an algebra~$A$ such that $\Conc A\cong S$. Given a variety~$\cV$ of algebras, \emph{the compact congruence class of~$\cV$}, denoted by $\Conc\cV$, is the class of all \jzs s isomorphic to $\Conc A$ for some $A\in\cV$. As illustrated by~\cite{PTW}, even the compact congruence classes of small varieties of lattices are complicated objects. For example, in case~$\cV$ is the variety of all lattices, $\Conc \cV$ contains all distributive \jzs s of cardinality at most~$\aleph_1$, but not all distributive \jzs s (cf.~\cite{CLP}).

Given varieties~$\cV$ and~$\cW$ of algebras, the \emph{critical point} of~$\cV$ and~$\cW$, denoted by $\crit{\cV}{\cW}$, is the smallest cardinality of a \jzs\ in $\Conc(\cV)-\Conc(\cW)$ if it exists, or~$\infty$, otherwise (i.e., if $\Conc\cV\subseteq\Conc\cW$).

Let~$I$ be a poset. A \emph{direct system indexed by~$I$} is a family $(A_i,f_{i,j})_{i\le j\text{ in }I}$ such that~$A_i$ is an algebra, $f_{i,j}\colon A_i\to A_j$ is a morphism of algebras, $f_{i,i}=\id_{A_i}$, and $f_{i,k}=f_{j,k}\circ f_{i,j}$, for all $i\le j\le k$ in~$I$.

Denote by~$\Sub V$ the subspace lattice of a vector space~$V$, and by~$\cM_n$ the variety of lattices generated by the lattice~$M_n$ of length two with~$n$ atoms, for $3\leq n\leq\omega$.
Using the theory of the \emph{dimension monoid} of a lattice, introduced by F. Wehrung in~\cite{Wehrung98}, together with some von Neumann regular ring theory, we prove in Section~\ref{S:Minor} that if~$\cV$ is a finitely generated variety of modular lattices with all subdirectly irreducible members of length less or equal to~$n$, then $\crit{\cV}{\Var(\Sub F^n)}\ge\aleph_2$ for any field~$F$. As an immediate application, $\crit{\cM_n}{\cM_3}\ge\aleph_2$ for every~$n$ with $3\leq n\leq\omega$ (cf. Corollary~\ref{C:Crit(MmMn)}). Thus, by using the result of M. Plo\v s\v cica in \cite{Ploscica00}, we obtain the equality $\crit{\cM_{m}}{\cM_n}=\aleph_2$ for all~$m$,~$n$ with $3\leq n<m\leq\omega$. Our proof does not rely on the approach used by Plo\v s\v cica in~\cite{Ploscica03} to prove the inequality $\crit{\cM_m^{0,1}}{\cM_n^{0,1}}\geq\aleph_2$, and it extends that result to the unbounded case. We also obtain a new proof of that result in Section~\ref{S:majoration-critpoint}, that does not even rely on the approach used by Plo\v s\v cica in~\cite{Ploscica00} to prove the inequality $\crit{\cM_{m}}{\cM_n}\leq\aleph_2$.

Let~$\cV$ be a variety of lattices, let~$\vec D$ be a diagram of \jzs s and \jzh s. A \emph{congruence-lifting} of~$\vec D$ in~$\cV$ is a diagram~$\vec L$ of~$\cV$ such that the composite $\Conc\circ\vec L$ is naturally equivalent to~$\vec D$.

In Section~\ref{S:majoration-critpoint}, we give a diagram of finite \jzs s that is congruence-liftable in~$\cM_n$, but not congruence-liftable in $\Var(\Sub F^3)$, for any finite field~$F$ and any~$n$ such that $2+\card F\le n\le\omega$. As the diagram of \jzs s is indexed by some ``good'' lattice, we obtain, using results of \cite{G:critpoint}, that $\crit{\cM_n}{\Var(\Sub F^3)}=\aleph_2$. This implies immediately that $\crit{\cM_4}{\cM_{3,3}}=\aleph_2$. Let~$F$ and~$K$ be finite fields such that $\card F>\card K$, we also obtain $\crit{\Var(\Sub F^3)}{\Var(\Sub K^3)}=\aleph_2$.

In a similar way, we prove that $\crit{\cM_\omega}{\cV}=\aleph_2$, for every finitely generated variety of lattices~$\cV$ such that $M_3\in\cV$.

\section{Basic concepts}\label{S:Basic}
We denote by $\dom f$ the domain of any function~$f$. A \emph{poset} is a partially ordered set.
Given a poset~$P$, we put
\[Q\dnw X=\setm{p\in Q}{(\exists x\in X)(p\le x)},\qquad Q\upw X=\setm {p\in Q}{(\exists x\in X)(p\ge x)},\]
for any $X, Q\subseteq P$, and we will write $\dnw X$ (resp., $\upw X$) instead of $P\dnw X$ (resp., $P\upw X$) in case~$P$ is understood. We shall also write $\dnw p$ instead of $\dnw\set{p}$, and so on, for $p\in P$. A poset~$P$ is \emph{lower finite} if $P\dnw p$ is finite for all $p\in P$. For $p, q\in P$ let $p\prec q$ hold, if $p<q$ and there is no $r\in P$ with $p<r<q$, in this case $p$ is called a \emph{lower cover} of~$q$. We denote by $P^=$ the set of all non-maximal elements in a poset~$P$. We denote by $\M(L)$ the set of all completely meet-irreducible elements of a lattice~$L$.

A \emph{$2$-ladder} is a lower finite lattice in which every element has at most two lower covers. S. Z. Ditor constructs in \cite{Ditor} a $2$-ladder of cardinality~$\aleph_1$.

For a set $X$ and a cardinal $\kappa$, we denote by:
\begin{align*}
[X]^\kappa & =\setm{Y\subseteq X}{\card Y=\kappa},\\
[X]^{\le\kappa} & =\setm{Y\subseteq X}{\card Y\le\kappa},\\
[X]^{<\kappa} & =\setm{Y\subseteq X}{\card Y<\kappa}.
\end{align*}

Denote by $\sP$ the category with objects the ordered pairs $(G,u)$ where $G$ is a pre-ordered abelian group and $u$ is an order-unit of $G$ (i.e., for each $x\in G$, there exists an integer~$n$ with $-nu\le x\le nu$), and morphisms $f\colon (G,u)\to (H,v)$ where $f\colon G\to H$ is an order-preserving group homomorphism and $f(u)=v$.

We denote by $\Dim$ the functor that maps a lattice to its \emph{dimension monoid}, introduced by F. Wehrung in \cite{Wehrung98}, we also denote by $\Delta(a,b)$ for $a\le b$ in~$L$ the canonical generators of $\Dim L$. We denote by $\GDim$ the functor that maps a lattice to the pre-ordered abelian universal group (also called Grothendieck group) of its dimension monoid. If~$L$ is a bounded lattice then (the canonical image in $\GDim(L)$ of) $\Delta(0_L,1_L)$ is an order-unit of $\GDim(L)$. If $f\colon L\to L'$ is a $0,1$-preserving homomorphism of bounded lattices, then $\GDim(f)\colon(\GDim(L),\Delta(0_L,1_L))\to(\GDim(L'),\Delta(0_{L'},1_{L'}))$ preserves the order-unit.

All our rings are associative but not necessarily unital.
\begin{itemize}
\item We denote by $\rideal(R)$ the poset of principal right ideals of every regular ring~$R$. The results of Fryer and Halperin in \cite[Section~3.2]{FrHa56}, imply that, $\rideal(R)$ is a~$0$-lattice, and for any homomorphism $f\colon R\to S$ of regular rings, the map $\rideal(f)\colon\rideal(R)\to\rideal(S)$, $I\mapsto f(I)S$ is a~$0$-lattice homomorphism (cf. Micol's thesis~\cite[Theorem~1.4]{Micol} for the unital case). Hence $\rideal$ is a functor from the category of regular rings to the category of~$0$-lattices with~$0$-lattice homomorphisms.
\item We denote by $V$ the functor from the category of unital rings with morphisms preserving units to the category of commutative monoids, that maps a unital ring~$R$ to the commutative monoid of all isomorphism classes of finitely generated projective right~$R$-modules and any homomorphism $f\colon R\to S$ of unital rings to the monoid homomorphism $V(f)\colon V(R)\to V(S)$, $\sum_i e_iR\mapsto \sum_i f(e_i)S$.
\end{itemize}
We denote by $\Id R$ (resp., $\Idc R$) the lattice of all two-sided ideals (resp., finitely generated two-sided ideals) of any ring~$R$. We denote by~$\Sub E$ the subspace lattice of a vector space~$E$. We denote by $\M_n(F)$ the~$F$-algebra of $n\times n$ matrices with entries from~$F$, for every field~$F$ and every positive integer~$n$. A \emph{matricial~$F$-algebra} is an~$F$-algebra of the form $\M_{k_1}(F)\times\dots\times \M_{k_n}(F)$, for positive integers $k_1,\dots,k_n$.

For a finitely generated projective right module~$P$ over a unital ring~$R$, we denote by $[P]$ the corresponding element in $K_0(R)$, that is, the stable isomorphism class of~$P$. We refer to \cite[Section~15]{Goodearl} for the required notions about the $K_0$ functor.

A \emph{$K_0$-lifting} of a pre-ordered abelian group with order-unit $(G,u)$ is a regular ring~$R$ such that $(K_0(R),[R])\cong(G,u)$. A $K_0$-lifting of a diagram $\vec G\colon I\to\sP$ is a diagram $\vec R\colon I\to\sP$ such that $(K_0(-),[-])\circ \vec R\cong\vec G$.

We denote by $\mtol$ the functor that sends a monoid to it maximal semilattice quotient, that is, $\mtol(M)=M/{\asymp}$ where $\asymp$ is the smallest congruence of $M$ such that $M/{\asymp}$ is a semilattice. We denote by $\pgrtol$ the functor that maps a partially pre-ordered abelian group $G$ to $\mtol (G^+)$ where $G^+$ is the monoid of all positive elements of $G$.

We denote by $\Var(L)$ (resp., $\Var_0(L)$, resp., $\Var_{0,1}(L)$) the variety of lattices (resp., lattices with~$0$, resp., bounded lattices) generated by a lattice~$L$.

A lattice~$K$ is a \emph{congruence-preserving extension of} a lattice~$L$, if~$L$ is a sublattice of~$K$ and $\Conc i\colon \Con L\to\Con K$ is an isomorphism, where $i\colon L\to K$ is the inclusion map.

We denote by $M_n$ and $M_{n,m}$ the lattices represented in Figure~\ref{F:MnandMnm}, for $3\le m,n\le\omega$, and by~$\cM_n$ and~$\cM_{n,m}$, respectively, the lattice varieties that they generate. We also denote by~$\cM_n^0$ the variety of lattices with~$0$ generated by $M_n$, and so on.

\setlength{\unitlength}{0.8mm}
\begin{figure}[bottom,here,top]
\begin{picture}(50,70)(-5,-5)
\put(20,0){\line(1,1){20}}
\put(20,0){\line(-1,1){20}}
\put(20,0){\line(-1,2){10}}
\put(19,-5){0}
\put(-6,20){$a_1$}
\put(4,20){$a_2$}
\put(34,20){$a_n$}
\put(19,42){1}

\put(0,20){\line(1,1){20}}
\put(10,20){\line(1,2){10}}
\put(40,20){\line(-1,1){20}}
\put(20,20){\dots}

\put(20,0){\circle*{2}}
\put(0,20){\circle*{2}}
\put(10,20){\circle*{2}}
\put(40,20){\circle*{2}}
\put(20,40){\circle*{2}}
\end{picture}\quad
\begin{picture}(70,70)(-5,-5)
\put(20,0){\line(1,1){20}}
\put(20,0){\line(-1,1){20}}
\put(20,0){\line(-1,2){10}}

\put(0,20){\line(1,1){20}}
\put(10,20){\line(1,2){10}}
\put(40,20){\line(-1,1){20}}
\put(20,20){\dots}

\put(40,20){\line(-1,2){10}}
\put(40,20){\line(1,1){20}}
\put(20,40){\line(1,1){20}}
\put(30,40){\line(1,2){10}}
\put(60,40){\line(-1,1){20}}
\put(40,40){\dots}

\put(20,0){\circle*{2}}
\put(0,20){\circle*{2}}
\put(10,20){\circle*{2}}
\put(40,20){\circle*{2}}
\put(20,40){\circle*{2}}
\put(30,40){\circle*{2}}
\put(60,40){\circle*{2}}
\put(40,60){\circle*{2}}
\put(19,-5){0}
\put(-6,20){$a_1$}
\put(4,20){$a_2$}
\put(34,20){$a_n$}
\put(14,40){$b_1$}
\put(24,40){$b_2$}
\put(54,40){$b_m$}
\put(38,62){1}
\end{picture}\caption{The lattices $M_n$ and $M_{n,m}$.}
\label{F:MnandMnm}
\end{figure}
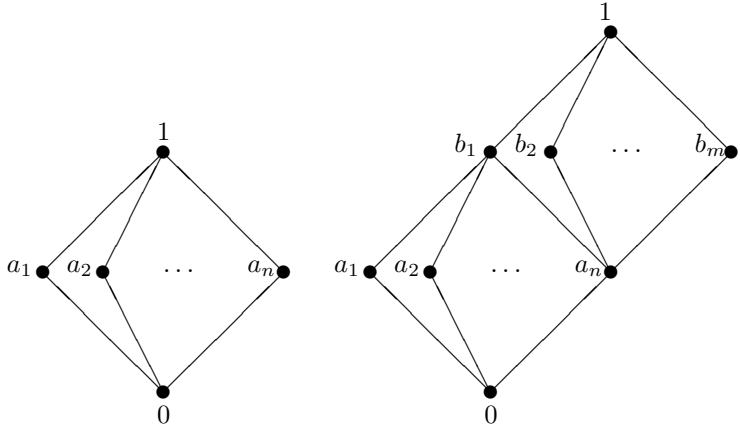

A lattice~$L$ \emph{satisfies Whitman's condition} if for all $a,b,c$, and $d$ in~$L$:
\[
a\wedge b\le c\vee d \text{\quad implies either $a\le c\vee d$\quad or\quad $b\le c\vee d$\quad or\quad $a\wedge b\le c$\quad or\quad $a\wedge b\le d$.}
\]
The lattice $M_n$ satisfies Whitman's condition for all $n\ge 3$.

\section{Lower bounds for some critical points}\label{S:Minor}

The following proposition is proved in \cite[Proposition~5.5]{Wehrung98}.
\begin{proposition}\label{P:W9855}
Let~$L$ be a modular lattice without infinite bounded chains. Let~$P$ be the set of all projectivity classes of prime intervals of~$L$. Given $\xi\in P$, denote by $\module{a,b}_\xi$ the number of occurrences of an interval in~$\xi$ in any maximal chain of the interval $[a,b]$. Then there exists an isomorphism $\pi\colon\Dim L\to (\mbZ^+)^{(P)}$ such that $\pi(\Delta(a,b))=\famm{\module{a,b}_\xi}{\xi\in P}$ for all $a\le b$ in~$L$.
\end{proposition}

This makes it possible to prove the following lemma, which gives an explicit description of $\GDim(L)$ for every modular lattice~$L$ of finite length (in such a case the set~$P$ is finite).

\begin{lemma}\label{L:GDimL-is-Zn}
Let~$L$ be a modular lattice of finite length, set $X=\M(\Con L)$. Then there exists an isomorphism
$\pi'\colon\GDim(L)\to \mbZ^X$ such that
\[\pi'(\Delta(a,b))=\famm{\lh([a/\theta,b/\theta])}{\theta\in X},\quad\text{for all $a\le b$ in~$L$}.\]
In particular $(\GDim(L),\Delta(0,1))$ is isomorphic to $(\mbZ^X,(\lh(L/\theta))_{\theta\in X})$.
\end{lemma}

\begin{proof}
Denote by~$P$ be the set of all projectivity classes of prime intervals of~$L$. For any $\xi\in P$ denote by $\theta_\xi$ the largest congruence of~$L$ that does not collapse any prime intervals in~$\xi$. As~$L$ is modular of finite length, the congruences of~$L$ are in one-to-one correspondence with subsets of $P$ (cf. \cite[Chapter~III]{GLT2}), and so the assignment $\xi\mapsto\theta_\xi$ defines a bijection from~$P$ onto~$X$. Moreover any prime interval not in~$\xi$ is collapsed by $\theta_\xi$, for any $\xi\in P$.
Let $a\le b$ in~$L$, let $\xi\in P$. Let $a_0\prec a_1\prec\dots\prec a_n$ in~$L$ such that $a_0=a$ and $a_n=b$. Let $0\le r_1<r_2<\dots<r_s<n$ be all the integers such that $[a_{r_k},a_{r_k+1}]\in\xi$ for all $1\le k\le s$. Thus $\module{a,b}_\xi=s$. Set $r_{s+1}=n$. As $[a_{r_k},a_{r_k+1}]\in\xi$ and $[a_{r_k+t},a_{r_k+t+1}]\not\in\xi$ for all $1\le t\le r_{k+1}-r_k-1$, we obtain that
\[a_{r_k}/\theta_\xi\prec a_{r_k+1}/\theta_\xi = a_{r_k+2}/\theta_\xi=\dots=a_{r_{k+1}}/\theta_\xi,\quad\text{for all $1\le k\le s.$}\]
Thus the following covering relations hold:
\[a/\theta_\xi=a_{r_1}/\theta_\xi\prec a_{r_2}/\theta_\xi\prec\dots\prec a_{r_s}/\theta_\xi\prec a_{r_{s+1}}/\theta_\xi=b/\theta_\xi.\]
So $\lh([a/\theta_\xi,b/\theta_\xi])=s=\module{a,b}_\xi$.
We conclude the proof by using Proposition~\ref{P:W9855}.
\end{proof}

\begin{proposition}\label{P:natequiv}
The following natural equivalences hold
\begin{align*}
(i)  & & \mtol\circ\Dim &\cong\Conc  & & \text{on lattices}\\
(ii) & & \mtol\circ V &\cong \Conc\circ \rideal  & & \text{on regular rings}
\end{align*}
\end{proposition}

\begin{proof}
$(i)$ follows from \cite[Corollary~2.3]{Wehrung98}, while $(ii)$ is contained in \cite[Corollary~2.23]{Goodearl}; see also the proof of \cite[Proposition~4.6]{Wehrung99}.
\end{proof}

We shall always apply this result to unital regular rings~$R$ such that $V(R)$ is cancellative (i.e.,~$R$ is unit-regular), so $K_0(R)^+=V(R)$, and to lattices~$L$ such that $\Dim L$ is cancellative, so $\GDim(L)^+\cong\Dim L$. Here $G^+$ denotes the positive cone of $G$, for any partially pre-ordered abelian group $G$.

The following theorem is proved in \cite[Theorem~15.23]{Goodearl}.
\begin{theorem}\label{T:goodearl}
Let~$F$ be a field, let~$R$ be a matricial~$F$-algebra, and let $S$ be a unit-regular~$F$-algebra.
\begin{enumerate}
\item Given any morphism $f\colon (K_0(R),[R])\to (K_0(S),[S])$ in $\sP$, the category of pre-ordered abelian groups with order-unit \textup(cf. Section~\ref{S:Basic}\textup), there exists an~$F$-algebra homomorphism $\phi\colon R\to S$ such that $K_0(\phi)=f$.
\item If $\phi,\psi\colon R\to S$ are~$F$-algebra homomorphisms, then $K_0(\phi)=K_0(\psi)$ if and only if there exists an inner automorphism $\theta$ of $S$ such that $\phi=\theta\circ\psi$.
\end{enumerate}
\end{theorem}

The following lemma is folklore.

\begin{lemma}\label{L:caract-K0-mat}
Let~$F$ be a field, let $\bu=(u_k)_{1\le k\le n}$ be a family of positive integers, let $R=\prod_{k=1}^n \M_{u_k}(F)$. Then $(K_0(R),[R])\cong(\mbZ^n,\bu)$.
\end{lemma}

\begin{lemma}\label{L:2-lader-of-Zn-is-liftable}
Let~$F$ be a field. Let~$I$ be a $2$-ladder, let $G_i=(\mbZ^{n_i},\bu^i=(u^i_k)_{1\le k\le n_i})$ such that $\bu^i$ is an order-unit, let $R_i=\prod_{k=1}^{n_i}\M_{u_k^i}(F)$ for all $i\in I$. Let $f_{i,j}\colon G_i\to G_j$ for all $i\le j$ in~$I$ such that $\vec{G}=(G_i,f_{i,j})_{i\le j\text{ in } I}$ is a direct system in $\sP$. Then there exists a direct system $(R_i,\phi_{i,j})_{i\le j\text{ in } I}$ of matricial~$F$-algebra which is a $K_0$-lifting of $(G_i,f_{i,j})_{i\le j\text{ in } I}$.
\end{lemma}

\begin{proof}
By Lemma~\ref{L:caract-K0-mat} there exists an isomorphism $\tau_i\colon(K_0(R_i),[R_i])\to G_i=(\mbZ^{n_i},\bu^i)$ in $\sP$, for all $i\in I$. Let $g_{i,j}=\tau_j^{-1}\circ f_{i,j}\circ \tau_i$, for all $i\le j$ in~$I$.

For $i=j=0$ (the smallest element of~$I$), we put $\phi_{0,0}=\id_{R_0}$. Let $i\in I$ with a lower cover $i'$. It follows from Theorem~\ref{T:goodearl}(1) that there exists $\psi_{i',i}\colon R_{i'}\to R_i$ such that $K_0(\psi_{i',i})=g_{i',i}$.

If $i$ has only $i'$ as lower cover, assume that we have a direct system $(R_j,\phi_{j,k})_{j\le k\le i'}$ lifting $(G_j,f_{j,k})_{j\le k\le i'}$. Set $\phi_{j,i}=\psi_{i',i}\circ \phi_{j,i'}$ for all $j<i$, and $\phi_{i,i}=\id_{R_i}$. It is easy to see that $(R_i,\phi_{j,k})_{j\le k\le i}$ is a direct system lifting $(G_j,f_{j,k})_{j\le k\le i}$.

Let $i$ has two distinct lower covers $i'$ and $i''$, and set $\ell=i'\wedge i''$. Assume that we have direct system $(R_j,\phi_{j,k})_{j\le k\le i'}$ and $(R_j,\phi_{j,k})_{j\le k\le i''}$ lifting $(G_j,f_{j,k})_{j\le k\le i'}$ and $(G_j,f_{j,k})_{j\le k\le i''}$ respectively. The following equalities hold
\[
K_0(\psi_{i',i}\circ \phi_{\ell,i'}) =K_0(\psi_{i',i})\circ K_0(\phi_{\ell,i'})= g_{i',i}\circ g_{\ell,i'}=g_{\ell,i}
\]

Similarly $K_0(\psi_{i'',i}\circ \phi_{\ell,i''})=g_{\ell,i}=K_0(\psi_{i',i}\circ \phi_{\ell,i'})$, thus, by Theorem~\ref{T:goodearl}(2), there exists an inner automorphism $\theta$ of $R_i$ such that $\theta\circ \psi_{i'',i}\circ \phi_{\ell,i''} = \psi_{i',i}\circ \phi_{\ell,i'}$. Put $\phi_{i',i} = \psi_{i',i}$ and $\phi_{i'',i}=\theta\circ \psi_{i'',i}$. Thus $\phi_{i',i}\circ\phi_{i'\wedge i'',i'}=\phi_{i'',i}\circ\phi_{i'\wedge i'',i''}$, so we can construct a direct system $(R_j,\phi_{j,k})_{j\le k\le i}$.

Hence, by induction, we obtain a direct system $(R_i,\phi_{i,j})_{i\le j\text{ in }I}$ of matricial~$F$-algebras, such that $K_0(\phi_{i,j})=g_{i,j}$ for all $i\le j$ in~$I$ as required.
\end{proof}

\begin{lemma}\label{L:lifting-with-regular-rings}
Let~$F$ be a field. Let~$L$ be a bounded modular lattice such that all finitely generated sublattices of~$L$ have finite length. Assume that $\card L\le\aleph_1$.
Then there exists a locally matricial ring~$R$ such that $\Con L\cong \Con \rideal(R)$ and $\rideal(R)\in\Var_{0,1}\famm{\Sub F^n}{n<\omega}$.

Moreover if there exists $n<\omega$ such that $n\ge \lh(K)$ for each simple lattice $K\in\Var(L)$ of finite length, then there exists a locally matricial ring~$R$ such that $\Con L\cong \Con \rideal(R)$ and $\rideal(R)\in\Var_{0,1}(\Sub F^n)$.
\end{lemma}

\begin{proof}
Let~$I$ be a $2$-ladder of cardinality~$\aleph_1$. Pick a surjection $\rho\colon I\tosurj L$ and denote by $L_i$ the sublattice of~$L$ generated by $\rho(I\dnw i)\cup\set{0,1}$, for each $i\in I$. Furthermore, denote by $f_{i,j}\colon L_i\to L_j$ the inclusion map, for all $i\le j$ in~$I$. Then $\vec L=(L_i,f_{i,j})_{i\le j\text{ in } I}$ is a direct system of modular lattices of finite length and $0,1$-lattice embeddings.

Assume that there exists $n<\omega$ such that $n\ge \lh(K)$ for each simple lattice $K\in\Var(L)$ of finite length.
Let $\vec G = \GDim\circ\vec{L}$, set $X_i=\M(\Con L_i)$ for all $i\in I$, and set $r^i_x=\lh(L_i/x)$ for each $x\in X_i$. The congruence lattice of any modular lattice of finite length is Boolean (cf. \cite[Chapter~III]{GLT2}), in particular, every subdirectly irreducible modular lattice of finite length is simple. This applies to the subdirectly irreducible lattice $L_i/x$, which is therefore simple. Thus $r^i_x \le n$, for all $i\in I$ and all $x\in X_i$. By Lemma~\ref{L:GDimL-is-Zn}, $G_i\cong (\mbZ^{X_i},(r^i_x)_{x\in X_i})$ for all $i\in I$.

Set $R_i=\prod_{x\in X_i}\M_{r^i_x}(F)$. By Lemma~\ref{L:caract-K0-mat}, $(K_0(R_i),[R_i])\cong (\mbZ^{X_i},(r^i_x)_{x\in X})\cong G_i$. By Lemma~\ref{L:2-lader-of-Zn-is-liftable}, there exists a direct system $\vec R=(R_i,\phi_{i,j})_{i\le j\text{ in }I}$ with morphisms preserving units, such that:
\begin{equation}\label{E:lifteq}
K_0\circ\vec R\cong\vec G = \GDim\circ\vec{L}.
\end{equation}
Moreover:
\[\rideal(R_i)\cong \rideal\left(\prod_{x\in X_i}\M_{r^i_x}(F)\right)\cong \prod_{x\in X_i}\rideal(\M_{r^i_x}(F))\cong \prod_{x\in X_i}\Sub F^{r^i_x}\in\Var_{0,1}(\Sub F^n).\]
Let $R=\varinjlim \vec R$. As $\rideal$ preserves direct limits, $\rideal(R)\cong \varinjlim(\rideal\circ\vec R)$, but $\rideal\circ\vec R$ is a diagram of $\Var_{0,1}(\Sub F^n)$, so $\rideal(R)\in\Var_{0,1}(\Sub F^n)$. Moreover the following isomorphisms hold:
\begin{align*}
\Conc\rideal(R) &\cong\pgrtol(K_0(R)) & &\text{by Proposition~\ref{P:natequiv}}\\
 &\cong\pgrtol(K_0(\varinjlim\vec R)) \\
&\cong\pgrtol(\varinjlim (K_0\circ \vec R)) & &\text{as $K_0$ preserves direct limits}\\
&\cong\pgrtol(\varinjlim (\GDim\circ \vec L)) & &\text{by \eqref{E:lifteq}}\\
&\cong\pgrtol(\GDim(\varinjlim\vec L)) & &\text{as $\GDim$ preserves direct limits}\\
&\cong\pgrtol(\GDim(L))\\
&\cong\Conc L & &\text{by Proposition~\ref{P:natequiv}.}
\end{align*}
The other case, without restriction on finite lengths of simple lattices, is similar.
\end{proof}

Lemma~\ref{L:lifting-with-regular-rings} works for bounded lattices, however any lattice can be embedded into a bounded lattice. In the rest of this section, using this result, we extend Lemma~\ref{L:lifting-with-regular-rings} to unbounded lattices.

\begin{lemma}\label{L:lattice-to-bounded-lattice}
Let~$L$ be a lattice, let $L'=L\sqcup\set{0,1}$ such that~$0$ is the smallest element of $L'$ and~$1$ is the largest. Let $f\colon L\toinj L'$ be the inclusion map. Then $\Conc f$ is a injective \jzh\ and $(\Conc f)(\Conc L)$ is an ideal of $\Conc L'$.
\end{lemma}

\begin{proof}
Let $\theta\in\Conc L$, let $L'_\theta=(L/\theta)\sqcup\set{0,1}$ such that~$0$ is the smallest element of~$L'_\theta$ and~$1$ is its largest element. The following map
\begin{align*}
g\colon L'&\to L'_\theta\\
x &\mapsto \begin{cases}
0 &\text{if $x=0$}\\
1 &\text{if $x=1$}\\
x/\theta &\text{if $x\in L$}
\end{cases}
\end{align*}
is a lattice homomorphism, and $\ker g=\theta\cup\set{(0,0),(1,1)}$, so the latter is a congruence of $L'$. It follows that $(\Conc f)(\theta)=\theta\cup\set{(0,0),(1,1)}$. Thus $\Conc f$ is an embedding. Let $\beta=\bigvee_{i=1}^n\Theta_{L'}(x_i,y_i)\in\Conc L'$, such that $\beta\subseteq(\Conc f)(\theta)$. We can assume that $x_i\not=y_i$ for all $1\le i\le n$. Thus, as $(x_i,y_i)\in\theta\cup\set{(0,0),(1,1)}$, $(x_i,y_i)\in\theta$ for all $1\le i\le n$. Let $\alpha=\bigvee_{i=1}^n\Theta_{L}(x_i,y_i)$, then $(\Conc f)(\alpha)=\beta$. Thus $(\Conc f)(\Conc L)$ is an ideal of $\Conc L'$.
\end{proof}

F. Wehrung proves the following proposition in \cite[Corollary~4.4]{Wehrung99}; the result also applies to the non-unital case, with a similar proof.

\begin{proposition}\label{P:weh99}
For any regular ring~$R$, $\Conc\rideal(R)$ is isomorphic to $\Idc R$.
\end{proposition}

\begin{lemma}
Let~$R$ be a regular ring, and let~$I$ be a two-sided ideal of~$R$. Then the following assertions hold
\begin{enumerate}
\item The set~$I$ is a regular subring of~$R$.
\item Any right \textup(resp., left\textup) ideal of~$I$ is a right \textup(resp., left\textup) ideal of~$R$.
\item In particular $\Id(I)=\Id(R)\dnw I$, and $\rideal(I)=\rideal(R)\dnw I$.
\end{enumerate}
\end{lemma}

\begin{proof}
The assertion $(1)$ follows from \cite[Lemma~1.3]{Goodearl}.

Let $J$ be a right ideal of~$I$, let $a\in J$, let $x\in R$. As~$I$ is regular there exists $y\in I$ such that $a=aya$, so $ax=ayax$, but $a\in I$, so $yax\in I$, moreover $J$ is a right ideal of~$I$, so $ax=ayax\in J$. Thus $J$ is a right ideal of~$R$. Similarly any left ideal of~$I$ is a left ideal of~$R$. Thus $\Id(I)=\Id(R)\dnw I$.

Let $a\in R$ idempotent. If $aR\subseteq I$, then $a\in I$, so $aI\subseteq aR = aaR\subseteq aI$, and so $aI=aR$, thus $aR\in\rideal(I)$. So $\rideal(I)=\rideal(R)\dnw I$.
\end{proof}

\begin{theorem}\label{T:crit-ge-aleph2}
Let~$F$ be a field. Let~$\cV$ be a variety of modular lattices \textup(resp., a variety of bounded modular lattices\textup). Assume that all finitely generated lattices of~$\cV$ have finite length. Then
\[\crit{\cV}{\Var_0\famm{\Sub F^n}{n\in\omega}}\ge\aleph_2 \quad\text{\textup(resp., }\crit{\cV}{\Var_{0,1}\famm{\Sub F^n}{n\in\omega}}\ge\aleph_2).
\]
Moreover for $L\in\cV$ of cardinality at most~$\aleph_1$, there exists a regular ring~$A$ such that $\Con L\cong \Con\rideal(A)$ and $\rideal(A)\in\Var_0\famm{\Sub F^n}{n\in\omega}$ \textup(resp., $\rideal(A)\in\Var_{0,1}\famm{\Sub F^n}{n\in\omega}$\textup).

If there exists $n<\omega$ such that $\lh(K)\le n$ for each simple lattice $K\in\cV$ of finite length, then:
\[\crit{\cV}{\Var_0(\Sub F^n)}\ge\aleph_2 \quad\text{\textup(resp., }\crit{\cV}{\Var_{0,1}(\Sub F^n)}\ge\aleph_2).
\]
Moreover for $L\in\cV$ of cardinality at most~$\aleph_1$, there exists a regular ring~$A$ such that $\Con L\cong \Con\rideal(A)$ and $\rideal(A)\in\Var_0(\Sub F^n)$ \textup(resp., $\rideal(A)\in\Var_{0,1}(\Sub F^n)$\textup).
\end{theorem}

Observe that $\rideal(A)$ is, in addition, relatively complemented; in particular, it is congruence-permutable.

\begin{proof}
The bounded case is an immediate application of Lemma~\ref{L:lifting-with-regular-rings}.

Let~$\cV$ be a variety of modular lattices in which finitely generated lattices have finite length. Let $L\in\cV$ such that $\card L\le\aleph_1$, let $L'=L\sqcup\set{0,1}$ as in Lemma~\ref{L:lattice-to-bounded-lattice} and let~$D$ be the ideal of $\Conc L'$ corresponding to $\Conc L$. By Chapter~I, Section~4, Exercise~14 in \cite{GLT2} we have $L'\in\cV$, thus, by Lemma~\ref{L:lifting-with-regular-rings}, there exists a regular ring~$R$ such that $\rideal(R)\in\Var_0(\Sub F^n)$, and $\Conc\rideal(R)\cong\Conc L'$. By Proposition~\ref{P:weh99}, $\Conc\rideal(R)\cong\Idc R$. Let~$I$ be the ideal of~$R$ corresponding to~$D$. Then $\Con L\cong\Id D\cong \Id R\dnw I\cong \Id I\cong \Con\rideal(I)$. Moreover $\rideal(I)= \rideal(R)\dnw I$ belongs to~$\cW$.
\end{proof}

We obtain the following generalization of M. Plo\v s\v cica's results in \cite{Ploscica03}.

\begin{corollary}\label{C:Crit(MmMn)}
Let~$m$,~$n$ be ordinals such that $3\leq n<m\leq\omega$. Then the equality $\crit{\cM_m}{\cM_n}=\aleph_2$ holds.
\end{corollary}

\begin{proof}
Every simple lattice of~$\cM_n$ has length at most two. Moreover, $\Sub \FF_2^2\cong M_3\in\cM_n$, where $\FF_2$ is the two-element field. Thus, by Theorem~\ref{T:crit-ge-aleph2}, $\crit{\cM_m}{\cM_n}\ge\aleph_2$.

Conversely, M. Plo\v s\v cica proves in \cite{Ploscica00} that there exists a \jz-semilattice of cardinality~$\aleph_2$, congruence-liftable in~$\cM_m$, but not congruence-liftable in~$\cM_n$. So $\crit{\cM_m}{\cM_n}\le\aleph_2$.
\end{proof}

In Section~\ref{S:majoration-critpoint} we shall give another \jzs\ of cardinality~$\aleph_2$, congruence-liftable in~$\cM_m$, but not congruence-liftable in~$\cM_n$.

\section{An upper bound of some critical points}\label{S:majoration-critpoint}

Using the results of \cite{G:critpoint}, we first prove that if a simple lattice of a variety of lattices~$\cV$ has larger length than all simple lattices of a finitely generated variety of lattices~$\cW$, then $\crit{\cV}{\cW}\le\aleph_0$.

\begin{remark}
Let $x\prec y$ in a lattice~$L$. Let $(\alpha_i)_{i\in I}$ be a family of congruences of~$L$, if $(x,y)\in\bigvee_{i\in I}\alpha_i$, then $(x,y)\in\alpha_i$ for some $i\in I$. In particular there exists a largest congruence separating $x$ and $y$. Such a congruence is completely \mirr, and in a lattice of finite height all completely \mirr\ congruences are of this form.
\end{remark}

\begin{lemma}\label{L:conddistributive}
Let~$L$ be a lattice and let $n\ge 0$. If $\Conc L\cong 2^n$ then $\lh(L)\ge n$. Moreover, if $C$ is a finite maximal chain of~$L$, then $\Conc f$ is surjective, where $f\colon C\to L$ is the inclusion map.
\end{lemma}

\begin{proof}
If~$L$ has no finite maximal chain then $\lh(L)\ge n$ is immediate. Assume that~$C$ is a finite maximal chain of~$L$. Denotes by $0=x_0\prec x_1\prec\dots\prec x_m=1$ the elements of~$C$. Denote by $f\colon C\to L$ the inclusion map.

Let $k\in\set{0,\dots,m-1}$. We have $x_k\prec x_{k+1}$, hence $\Theta_L(x_k,x_{k+1})$ is join-irreducible in $\Conc L$. As $\Conc L$ is Boolean, $\Theta_L(x_k,x_{k+1})$ is an atom of $\Conc L$.

Let $\theta$ be an atom of $\Conc L$, we have:
\[
\theta\le \Theta_L(0,1)=\bigvee_{k=0}^{m-1}\Theta_L(x_k,x_{k+1})
\]
So there exists $k\in\set{0,\dots,m-1}$ such that $\theta\le \Theta_L(x_k,x_{k+1})$. As $\Theta_L(x_k,x_{k+1})$ is an atom of $\Conc L$, we have $\theta=\Theta_L(x_k,x_{k+1})$. It follows that $\Conc f$ is surjective, so $m\ge n$ and so $\lh(L)\ge n$.
\end{proof}

\begin{theorem}\label{T:thealeph0}
Let~$\cV$ be a variety of lattices \textup(resp., a variety of bounded lattices\textup), let~$\cW$ be a finitely generated variety of lattices, let~$D$ be a finite \jzs. If there exists a lifting $K\in\cV$ of $D$ of length greater than every lifting of $D$ in~$\cW$, then $\crit{\cV}{\cW}\le\aleph_0$. Moreover if~$\cV$ is a finitely generated variety of modular lattices and~$\cW$ is not trivial, then $\crit{\cV}{\cW}=\aleph_0$.
\end{theorem}

\begin{proof}
As $D$ is finite, taking a sublattice, we can assume that $\card K\le\aleph_0$. Let~$n$ be the greatest length of a lifting of $D$ in~$\cW$. As $\lh(K)>n$, there exists a chain~$C$ of~$K$ of length $n+1$ (resp., we can assume that $C$ has~$0$ and~$1$). Let $f\colon C\to K$ be the inclusion map. Assume that there exists a lifting $g\colon C'\to K'$ of $\Conc f$ in~$\cW$. As $f$ is an embedding, $g$ is also an embedding. As $\Conc K'\cong\Conc K\cong D$, $\lh(K')\le n$. Moreover $\Conc C'\cong\Conc C\cong 2^{n+1}$, thus, by Lemma~\ref{L:conddistributive}, $\lh(C')=n+1$. So $n\ge\lh(K')\ge\lh(C')=n+1$; a contradiction.

Therefore $\Conc f$ has no lifting in~$\cW$. So, as $\card K\le\aleph_0$ and by \cite[Corollary~7.6]{G:critpoint}, $\crit{\cV}{\cW}\le\aleph_0$ (in the bounded case $f$ preserves bounds, thus the result of \cite{G:critpoint} also applies).

Moreover if~$\cV$ is a finitely generated variety of modular lattices, then the finite \jzs s with congruence-lifting in~$\cV$ are the finite Boolean lattices. Finite Boolean lattices are also liftable in~$\cW$. Hence $\crit{\cV}{\cW}=\aleph_0$.
\end{proof}

The following corollary is an immediate application of Theorem~\ref{T:thealeph0} and Theorem~\ref{T:crit-ge-aleph2}. It shows that the critical point between a finitely generated variety of modular lattices and a variety generated by a lattice of subspaces of a finite vector space, cannot be~$\aleph_1$.

\begin{corollary}
Let~$\cV$ be a finitely generated variety of modular lattices, let~$F$ be a finite field, let $n\ge 1$ be an integer. If there exists a simple lattice in $K\in\cV$ such that $\lh(K)>n$, then $\crit{\cV}{\Var(\Sub F^n)}=\aleph_0$, else $\crit{\cV}{\Var(\Sub F^n)}\ge\aleph_2$.
\end{corollary}

We shall now give a diagram of \jzs s $\vec S$, congruence-liftable in~$\cM_n$, such that for every finitely generated variety~$\cV$, generated by lattices of length at most three, the diagram $\vec S$ is congruence-liftable in~$\cV$ if and only if $M_n\in\cV$.

Let $n\ge 3$ be an integer. Set $\un=\set{0,1,\dots,n-1}$, and set:
\[
I_n=\setm{P\in \Pow(\un)}{\text{either $\card(P)\le 2$ or $P=\un$}}.
\]

Denote by $a_0,\dots,a_{n-1}$ the atoms of $M_n$. Set $A_P=\setm{a_x}{x\in P}\cup\set{0,1}$, for all $P\in I_n$. Let $f_{P,Q}\colon A_P\to A_Q$ be the inclusion map for all $P\subseteq Q$ in~$I_n$. Then $\vec A=(A_P,f_{P,Q})_{P\subseteq Q\text{ in }I_n}$ is a direct system in~$\cM_n^{0,1}$. The diagram~$\vec S $ is defined as $\Conc\circ\vec A$.

\begin{lemma}\label{L:cong-for-Mn}
Let $\vec B=(B_P,g_{P,Q})_{P\subseteq Q\text{ in }I_n}$ be a congruence-lifting of $\Conc\circ \vec A$ by lattices, with all the maps $g_{P,Q}$ inclusion maps, for all $P\subseteq Q$ in $I_n$. Let $u<v$ in $B_\emptyset$. Let $P\in I_n$ then:
\[
\Theta_{B_P}(u,v)=B_P\times B_P,\quad\text{the largest congruence of $B_P$}.
\]
Let $\vec\xi=(\xi_P)_{P\in I_n}\colon \Conc\circ\vec A\to \Conc\circ\vec B$ be a natural equivalence. Let $x,y\in\un$ distinct. Let $b_x\in [u,v]_{B_{\set{x}}}$ and $b_y\in [u,v]_{B_{\set{y}}}$. Set $P=\set{x,y}$. Let $c\in\set{0,1}$. Then the following assertions hold:
\begin{enumerate}
\item If $\Theta_{B_{\set{x}}}(u,b_x) = \xi_{\set{x}}(\Theta_{A_{\set{x}}}(c,a_x))$, then $\Theta_{B_P}(u,b_x) = \xi_P(\Theta_{A_P}(c,a_x))$.
\item If $\Theta_{B_{\set{z}}}(u,b_z) = \xi_{\set{z}}(\Theta_{A_{\set{z}}}(c,a_z))$ for all $z\in\set{x,y}$, then $b_x\wedge b_y=u$.
\item If $\Theta_{B_{\set{x}}}(b_x,v) = \xi_{\set{x}}(\Theta_{A_{\set{x}}}(c,a_x))$, then $\Theta_{B_P}(b_x,v) = \xi_P(\Theta_{A_P}(c,a_x))$.
\item If $\Theta_{B_{\set{z}}}(b_z,v) = \xi_{\set{z}}(\Theta_{A_{\set{z}}}(c,a_z))$ for all $z\in\set{x,y}$, then $b_x\vee b_y=v$.
\item If $\Theta_{B_{\set{x}}}(u,b_x) = \xi_{\set{x}}(\Theta_{A_{\set{x}}}(c,a_x))$ and $\Theta_{B_{\set{y}}}(b_y,v) = \xi_{\set{y}}(\Theta_{A_{\set{y}}}(c,a_y))$, then $b_x\le b_y$.
\end{enumerate}
\end{lemma}

\begin{proof}
As $f_{P,Q}$ preserves bounds, $\Conc f_{P,Q}$ preserves bounds, thus $\Conc g_{P,Q}$ preserves bounds, for all $P\subseteq Q$ in $I_n$. Let $u<v$ in $B_\emptyset$. As $B_\emptyset$ is simple, $\Theta_{B_\emptyset}(u,v)$ is the largest congruence of $B_\emptyset$. Moreover, $\Conc g_{\emptyset,P}$ preserves bounds, for all $P\in I_n$. Hence:
\[
\Theta_{B_P}(u,v)=B_P\times B_P,\quad\text{the largest congruence of $B_P$}.
\]
\begin{enumerate}
\item The following equalities hold:
\begin{align*}
\Theta_{B_P}(u,b_x) &= (\Conc g_{\set{x},P})(\Theta_{B_{\set{x}}} (u,b_x))\\
&= (\Conc g_{\set{x},P})(\xi_{\set{x}}(\Theta_{A_{\set{x}}}(c,a_x))) &\text{by assumption}\\
&= \xi_P \circ (\Conc f_{\set{x},P})(\Theta_{A_{\set{x}}}(c,a_x))\\
&= \xi_P(\Theta_{A_P}(c,a_x)).
\end{align*}
\item The following containments hold:
\begin{align*}
\Theta_{B_P}(u,b_x\wedge b_y) &\subseteq \Theta_{B_P}(u,b_x)\cap \Theta_{B_P}(u,b_y)\\
&= \xi_P(\Theta_{A_P}(c,a_x))\cap \xi_P(\Theta_{A_P}(c,a_y)) &\text{by $(1)$}\\
&= \xi_P(\Theta_{A_P}(c,a_x)\cap \Theta_{A_P}(c,a_y))\\
&= \xi_P(\id_{A_P})=\id_{B_P}.
\end{align*}
so $u=b_x\wedge b_y$.
\item Similar to $(1)$.
\item Similar to $(2)$.
\item The following containments hold:
\begin{align*}
\Theta_{B_P}(b_y,b_x\vee b_y) &\subseteq \Theta_{B_P}(u,b_x)\cap \Theta_{B_P}(b_y,v)\\
&= \xi_P(\Theta_{A_P}(c,a_x))\cap \xi_P(\Theta_{A_P}(c,a_y)) &\text{by $(1)$ and $(3)$}\\
&= \xi_P(\Theta_{A_P}(c,a_x)\cap \Theta_{A_P}(c,a_y))\\
&= \xi_P(\id_{A_P})=\id_{B_P}.
\end{align*}
so $b_y=b_x\vee b_y$, thus $b_x\le b_y$.\qed
\end{enumerate}
\renewcommand{\qed}{}
\end{proof}

The following lemma shows that if we have some ``small'' enough congruence-lifting of $\Conc\circ \vec A$ in a variety, then $M_n$ belongs to this variety.

\begin{lemma}\label{L:lifting-imply-Mn}
Let $\vec B=(B_P,g_{P,Q})_{P\subseteq Q\text{ in }I_n}$ be a congruence-lifting of $\Conc\circ \vec A$ by lattices. Assume that~$B_{\set{x}}$ is a chain of length two for all $x\in \un$. Then $M_n$ can be embedded into $B_{\un}$.
\end{lemma}

\begin{proof}
Let $\vec\xi=(\xi_P)_{P\in I_n}\colon\Conc\circ\vec A\to\Conc\circ\vec B$ be a natural equivalence. As $f_{P,Q}$ is an embedding, $\Conc f_{P,Q}$ separates~$0$, so $\Conc g_{P,Q}$ separates~$0$, hence $g_{P,Q}$ is an embedding, thus we can assume that $g_{P,Q}$ is the inclusion map from $B_P$ into $B_Q$, for all $P\subseteq Q$ in~$I_n$.

Let $u<v$ in $B_\emptyset$. By Lemma~\ref{L:cong-for-Mn}, $\Theta_{B_{\set{x}}}(u,v)$ is the largest congruence of $B_{\set{x}}$. Moreover $B_{\set{x}}$ is the $3$-element chain, so $u$ is the smallest element of $B_{\set{x}}$ while $v$ is its largest element. Denote by $b_x$ the middle element of $B_{\set{x}}$.

The congruence $\xi_{\set{x}}(\Theta_{A_{\set{x}}}(0,a_x))$ is join-irreducible, thus it is either equal to $\Theta_{B_{\set{x}}}(u,b_x)$ or to $\Theta_{B_{\set{x}}}(b_x,v)$. Set:
\[
X'=\setm{x\in \un}{\xi_{\set{x}}(\Theta_{A_{\set{x}}}(0,a_x)) = \Theta_{B_{\set{x}}}(u,b_x)},\]
\[X''=\setm{x\in \un}{\xi_{\set{x}}(\Theta_{A_{\set{x}}}(0,a_x)) = \Theta_{B_{\set{x}}}(b_x,v)}.\]
As $\Theta_{A_{\set{x}}}(0,a_x)$ is the complement of $\Theta_{A_{\set{x}}}(a_x,1)$ and $\Theta_{B_{\set{x}}}(u,b_x)$ is the complement of $\Theta_{B_{\set{x}}}(b_x,v)$, we also get that:
\[
X'=\setm{x\in \un}{\xi_{\set{x}}(\Theta_{A_{\set{x}}}(a_x,1)) = \Theta_{B_{\set{x}}}(b_x,v)}\]
\[X''=\setm{x\in \un}{\xi_{\set{x}}(\Theta_{A_{\set{x}}}(a_x,1)) = \Theta_{B_{\set{x}}}(u,b_x)}.\]
Moreover $\un=X'\cup X''$. As $\card \un\ge 3$, either $\card X'\ge 2$ or $\card X''\ge 2$.

Assume that $\card X'\ge 2$. Let $x,y$ in $X'$ distinct. By Lemma~\ref{L:cong-for-Mn}(2), $b_x\wedge b_y=u$. By Lemma~\ref{L:cong-for-Mn}(4), $b_x\vee b_y=v$.

Now assume that $X''\not=\emptyset$. Let $z\in X''$. As $\xi_{\set{x}}(\Theta_{A_{\set{x}}}(0,a_x))=\Theta_{B_{\set{x}}}(u,b_x)$ and $\xi_{\set{z}}(\Theta_{A_{\set{z}}}(0,a_z))=\Theta_{B_{\set{z}}}(b_z,v)$, it follows from Lemma~\ref{L:cong-for-Mn}(5) that $b_x\le b_z$. Similarly, as $\xi_{\set{z}}(\Theta_{A_{\set{z}}}(a_z,1))=\Theta_{B_{\set{z}}}(u,b_z)$ and $\xi_{\set{y}}(\Theta_{A_{\set{y}}}(a_y,1))=\Theta_{B_{\set{y}}}(b_y,v)$, it follows from Lemma~\ref{L:cong-for-Mn}(5) that $b_z\le b_y$. Thus $b_x\le b_y$. So $u=b_x\wedge b_y=b_x>u$, a contradiction.

Thus $X''=\emptyset$, so $X'=\un$, and so $\set{u,b_0,b_1,\dots,b_n,v}$ is a sublattice of $B_{\un}$ isomorphic to $M_n$. The case $\card X''\ge 2$ is similar.
\end{proof}

We shall now use a tool introduced in \cite{G:critpoint} to prove that having a congruence-lifting of $\Conc\circ\vec A$ is equivalent to having a congruence-lifting of some \jzs\ of cardinality~$\aleph_2$. This requires the following infinite combinatorial property, proved by A. Hajnal and A. M\'at\'e in \cite{HajMat}, see also \cite[Theorem~46.2]{EHMR}. This property is also used by M. Plo\v s\v cica in \cite{Ploscica00}.

\begin{proposition}\label{P:HajMat}
Let $n\ge 0$ be an integer, let $\alpha$ be an ordinal, let $\kappa\ge\aleph_{\alpha+2}$, let $f\colon[\kappa]^2\to[\kappa]^{<\aleph_\alpha}$. Then there exists $Y\in[\kappa]^n$ such that $a\not\in f(\set{b,c})$ for all distinct $a,b,c\in Y$.
\end{proposition}

Now recall the definition of supported poset and norm-covering introduced in \cite[Section~4]{G:critpoint}.
\begin{definition}\label{D:KerSupp}
A finite subset $V$ of a poset $U$ is a \emph{kernel}, if for every $u\in U$,
there exists a largest element $v\in V$ such that $v\le u$. We denote
this element by $V\cdot u$.

We say that $U$ is \emph{supported}, if every finite subset of $U$ is
contained in a kernel of~$U$.
\end{definition}

We denote by $V\cdot \bu$ the largest element of $V\cap\bu$, for every
kernel $V$ of $U$ and every ideal~$\bu$ of~$U$.
As an immediate application of the finiteness of kernels, we obtain
that any intersection of a nonempty set of kernels of a poset $U$ is a kernel of $U$.

\begin{definition}\label{D:normcovering}
A \emph{norm-covering} of a poset~$I$ is a pair $(U,\module{\cdot})$, where $U$ is a supported poset and $\module{\cdot}\colon U\to I$, $u\mapsto\module{u}$ is an order-preserving map.

A \emph{sharp ideal} of $(U,\module{\cdot})$ is an ideal $\bu$ of $U$ such
that $\setm{\module{v}}{v\in \bu}$ has a largest element. For example, for every $u\in U$, the principal ideal $U\dnw u$ is sharp; we shall often identify~$u$ and~$U\dnw u$. We denote
this element by $\module{\bu}$. We denote by $\Ids(U,\module{\cdot})$ the set
of all sharp ideals of $(U,\module{\cdot})$, partially ordered by inclusion.

A sharp ideal $\bu$ of $(U,\module{\cdot})$ is \emph{extreme},
if there is no sharp ideal $\bv$ with $\bv>\bu$ and $\module{\bv}=\module{\bu}$. We
denote by $\Ide(U,\module{\cdot})$ the set of all extreme ideals of $(U,\module{\cdot})$.

Let $\kappa$ be a cardinal number. We say that $(U,\module{\cdot})$ is \emph{$\kappa$-compatible},
if for every order-preserving map $F\colon\Ide(U,\module{\cdot})\to\mathfrak{P}(U)$ such that $\card F(\bu)<\kappa$ for all $\bu\in \Ide(U,\module{\cdot})^=$, there exists an order-preserving map $\sigma\colon I\to \Ide(U,\module{\cdot})$ such that:
\begin{enumerate}
\item The equality $\module{\sigma(i)}=i$ holds for all $i\in I$.
\item The containment $F(\sigma(i))\cap\sigma(j)\subseteq\sigma(i)$ holds for all $i\leq j$ in~$I$.
\end{enumerate}
\end{definition}

\begin{lemma}\label{L:supported}
Let $X$ be a set, let $(A_x)_{x\in X}$ be a family of sets, let:
\[
U=\bigsqcup_{P\in \Powf{X}}\prod_{x\in P} A_x.
\]
We view the elements of $U$ as \textup(partial\textup) functions and ``to be greater" means ``to extend". Then $U$ is a supported poset.
\end{lemma}

\begin{proof}
Let $V$ be a finite subset of $U$. Let $Y_x=\setm{u_x}{u\in V\text{ and }x\in\dom u}$ for all $x\in X$. Let $D=\bigcup_{u\in V}\dom u$. Let:
\[W=\setm{u\in U}{\dom u\subseteq D\text{ and }(\forall x\in\dom u)(u_x\in Y_x)}\]
the set $D$, and the sets $Y_x$ for $x\in X$ are all finite, so $W$ is finite.

Let $u\in U$, let $P=\setm{x\in\dom u}{x\in D\text{ and }u_x\in Y_x}$. Then $u\res P \in W$. Moreover let $w\in W$ such that $w\le u$. Let $x\in\dom w$, then $x\in D$, and $u_x=w_x\in Y_x$, thus $\dom w\subseteq P$, so $w\le u\res P$. Therefore $u\res P$ is the largest element of $W\dnw u$.
\end{proof}

Using Lemma~\ref{L:supported} and Proposition~\ref{P:HajMat} we can construct a~$\aleph_\alpha$-compatible lower finite norm-covering of~$I_n$, the poset constructed earlier.

\begin{lemma}\label{L:I-compatible-norm-covering}
Let $\alpha$ be an ordinal. Let $U=\bigsqcup_{P\in\Pow(\un)}\aleph_{\alpha+2}^P$, partially ordered by inclusion. Let
\begin{align*}
\module{\cdot}\colon U &\to I_n\\
u &\mapsto \module{u}=
\begin{cases}
\dom u &\text{if $\card (\dom u)\le 2$}\\
\un &\text{otherwise}.
\end{cases}
\end{align*}
Then $(U,\module{\cdot})$ is a~$\aleph_\alpha$-compatible lower finite norm-covering of~$I_n$. Moreover $\card U=\aleph_{\alpha+2}$.
\end{lemma}

\begin{proof}
By Lemma~\ref{L:supported}, the set $U$ is supported. Moreover $\module{\cdot}$ preserves order, so $(U,\module{\cdot})$ is a norm-covering of~$I_n$. The poset $U$ is lower finite.

Extreme ideals are of the form $\dnw u$, where $u\in U$ and $\dom u\in I_n$, so we identify the corresponding extreme ideal with $u$. Thus $\Ide(U,\module{\cdot})=\setm{u\in U}{\dom u\in I_n}$.

Let $F\colon \Ide(U,\module{\cdot})\to \Pow(U)$ be an order-preserving map such that $\card F(\bu)<\aleph_\alpha$ for all $\bu\in \Ide(U,\module{\cdot})^=$, let
\begin{align*}
G\colon[\aleph_{\alpha+2}]^2 &\to[{\aleph_{\alpha+2}}]^{<\aleph_\alpha}\\
s&\mapsto\bigcup\Setm{\im v}{u \in\bigcup_{P\in I_n-\set{\un}}s^P\text{ and }v\in F(u)}.
\end{align*}
By Proposition~\ref{P:HajMat}, there exists $A\subset\aleph_{\alpha+2}$ such that $\card{A}=n$ and $a\not\in G(\set{b,c})$ for all distinct $a,b,c\in A$. Let $u\colon \un\to A$ be a one-to-one map. Let $\phi\colon I_n\to \Ide(U,\module{\cdot})$, $P\mapsto u\res P$. Then $\module{\phi(P)}=P$. Let $P\subsetneq Q$ in~$I_n$, let $v\in F(u\res P)\dnw (u\res Q)$. Let $x\in \dom v-P$. As $P\in I_n$, and $P\not=\un$, $\card P\le 2$. Let $P'=\set{y,z}\subseteq \un$, such that $y,z$ are distinct, $P\subseteq P'$, and $x\not\in P'$. Let $s=\set{u_y,u_z}$, then $u\res P'\in s^{P'}$, as $v\in F(u\res P)\subseteq F(u\res P')$, $v_x\in G(s)$. Moreover $u_x,u_y,u_z\in A$ are distinct, thus $u_x\not\in G(\set{u_y,u_z})=G(s)$, so $v_x\not=u_x$ in contradiction with $v\le u$, so $\dom v\subseteq P$, and so $v\le u\res P$.
\end{proof}

Using the results of \cite{G:critpoint} together with Lemma~\ref{L:I-compatible-norm-covering}, we obtain the following result.

\begin{lemma}\label{L:lifting-diagram-idexed-by-I}
Let~$\cV$ be a variety of algebras with a countable similarity type, let~$\cW$ be a finitely generated congruence-distributive variety such that $\crit{\cV}{\cW}>\aleph_2$. Let $\vec{D}\colon I_n\to\cS$ be a diagram of finite \jzs s. If $\vec D$ is congruence-liftable in~$\cV$, then $\vec D$ is congruence-liftable in~$\cW$.
\end{lemma}

\begin{proof}
In this proof we use, but do not give, many definitions of \cite{G:critpoint}. By Lemma~\ref{L:I-compatible-norm-covering} there exists $(U,\module{\cdot})$ a~$\aleph_0$-compatible lower finite norm-covering of~$I_n$ such that $\card U=\aleph_{2}$. Let $J$ be a one-element ordered set. By \cite[Lemma~3.9]{G:critpoint},~$\cW$ is $(\Ide(U,\module{\cdot})^=,J,\aleph_0)$-L\"owenheim-Skolem.

Let $\vec A=(A_P,f_{P,Q})_{P\subseteq Q\text{ in }I_n}$ be a congruence-lifting of $\vec D$ in~$\cV$. As $\Conc A_P$ is finite, using \cite[Lemma~3.6]{G:critpoint}, taking sublattices we can assume that $A_P$ is countable for all $P\in I_n$. By~\cite[Lemma~6.7]{G:critpoint}, there exists an $U$-quasi-lifting $(\tau,\Cond(\vec A,U))$ of $\vec D$ in~$\cV$. Moreover:
\[
\card\Cond(\vec A,U) \le \sum_{V\in\Powf{U}}\card\left(\prod_{u \in V}A_{\module u}\right)\le \sum_{V\in\Powf{U}}\aleph_0\le\aleph_2
\]
As $\crit{\cV}{\cW}>\aleph_2$, there are $B\in\cW$ and an isomorphism $\xi\colon\Conc\Cond(\vec A,U)\to\Conc B$.
So $(\tau\circ\xi^{-1},B)$ is an $U$-quasi-lifting of $\vec D$. Moreover~$\cW$ is $(\Ide(U,\module{\cdot})^=,J,\aleph_0)$-L\"owenheim-Skolem, hence, by \cite[Theorem~6.9]{G:critpoint}, with $I=I_n$, there exists a congruence-lifting of~$\vec D$ in~$\cW$.
\end{proof}

A similar proof, using Lemma~3.6, Lemma~3.7, Lemma~6.7, and Theorem~6.9 in \cite{G:critpoint} together with Lemma~\ref{L:I-compatible-norm-covering}, yields the following generalization of Lemma~\ref{L:lifting-diagram-idexed-by-I}.

\begin{lemma}\label{L:lifting-diagram-idexed-by-I-nonFG}
Let $\alpha\ge 1$ be an ordinal. Let~$\cV$ and~$\cW$ be varieties of algebras, with similarity types of cardinality $<\aleph_\alpha$. Let $\vec{D}=(D_P,\varphi_{P,Q})_{P\subseteq Q\text{ in }I_n}$ be a direct system of \jzs s. Assume that the following conditions hold:
\begin{enumerate}
\item $\crit{\cV}{\cW}>\aleph_{\alpha+2}$.
\item $\card(D_P)< \aleph_\alpha$, for all $P\in I_n-\set{\un}$.
\item $\card(D_{\un})\le\aleph_{\alpha+2}$.
\item $\vec D$ is congruence-liftable in~$\cV$.
\end{enumerate}
Then $\vec D$ is congruence-liftable in~$\cW$.
\end{lemma}

The following lemma implies, in particular, that a modular lattice of length three is a congruence-preserving extension of one of its subchains.

\begin{lemma}\label{L:lattice-to-chain}
Let~$L$ be a lattice of length at most three, let $u,v$ in~$L$ such that $\Theta_L(u,v)=L\times L$. If $\Conc L\cong 2^2$, then there exists $x\in L$ with $u<x<v$ such that~$L$ is a congruence-preserving extension of the chain $C=\set{u,x,v}$.
\end{lemma}

\begin{figure}[here,top,bottom]
\setlength{\unitlength}{0.7mm}
\begin{picture}(50,50)(-5,-5)
\put(20,0){\line(-1,1){20}}
\put(20,0){\line(1,1){15}}
\put(35,15){\line(0,1){10}}
\put(0,20){\line(1,1){20}}
\put(35,25){\line(-1,1){15}}
\put(20,0){\circle*{2}}
\put(0,20){\circle*{2}}
\put(35,15){\circle*{2}}
\put(35,25){\circle*{2}}
\put(20,40){\circle*{2}}
\put(19,-4){$u$}
\put(19,42){$v$}
\put(-5,19){$x$}
\put(38,14){$y$}
\put(38,24){$z$}
\end{picture}\caption{The lattice $N_5$.}
\label{F:N5}
\end{figure}
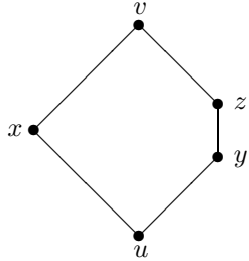

\begin{proof}
As $\Conc L\cong 2^2$, $\lh([u,v])\ge 2$. If $\lh([u,v])=2$, then let $C=\set{u,x,v}$, where $x$ is any element such that $u<x<v$. Let $i\colon C\to L$ the inclusion map. The morphism $\Conc i\colon \Conc C\to \Conc L$ is onto, moreover $\Conc C\cong 2^2\cong\Conc L$, so $\Conc i$ is an isomorphism.

Now assume that $[u,v]$ has length three. As $\lh(L)\le 3$, $\lh(L)=3$, $u$ is the smallest element of~$L$, and $v$ is the largest element.

Assume that~$L$ has a sublattice isomorphic to $N_5$, as labeled in Figure~\ref{F:N5}. Then $C=\set{u,y,z,v}$ is a maximal chain of~$L$. Let $i\colon C\to L$ be the inclusion map. By Lemma~\ref{L:conddistributive}, $\Conc i$ is surjective. Thus, as $\Con L\cong 2^2$, and $\Theta_L(u,y)$, $\Theta_L(y,z)$, and $\Theta_L(z,v)$ are all the atoms of $\Con L$,
\[\Theta_L(y,z)\subseteq\Theta_L(u,y)\cap\Theta_L(y,z)\cap\Theta_L(z,v) = \id_L,\]
a contradiction. Thus~$L$ does not contain any lattice isomorphic to $N_5$, that is,~$L$ is modular.

As $\Con L\cong 2^2$ and $\lh(L)=3$,~$L$ is not distributive. Hence there exists a sublattice of~$L$ isomorphic to $M_3$, let $a< x_1,x_2,x_3<b$ be its elements. As~$L$ is modular, $[a,x_1]_L\cong[x_1,b]_L$, thus $\lh([a,b]_L)$ is even. But $2\le \lh([a,b]_L)\le 3$, so $\lh([a,b]_L)=2$, thus $a\prec x_1\prec b$. This chain can be completed into a maximal chain $c\prec a\prec x_1\prec b$ or $a\prec x_1\prec b \prec c$. By symmetry, we may assume that $b<c$. Observe that $a=u$ and $c=v$. Set $C=\set{u,b,v}$ and $C_1=\set{u,x_1,b,v}$. Let $i\colon C\to L$ and $i_1\colon C_1\to L$ be the inclusion maps. As $C_1$ is a maximal chain, $\Conc i_1$ is onto. As $\Theta_L(u,x_1)=\Theta_L(x_1,b)=\Theta_L(u,b)$, $\Conc i_1$ and $\Conc i$ have the same image, thus $\Conc i$ is onto, so $\Conc i$ is an isomorphism.
\end{proof}

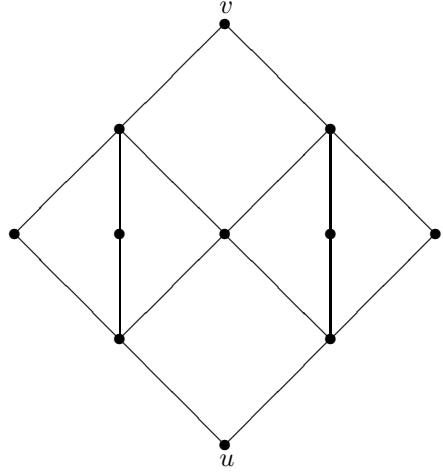
\begin{figure}[here,top,bottom]
\setlength{\unitlength}{0.7mm}
\begin{picture}(90,90)(-5,-5)
\put(40,0){\line(1,1){40}}
\put(40,0){\line(-1,1){40}}

\put(60,20){\line(0,1){40}}
\put(20,20){\line(0,1){40}}

\put(20,20){\line(1,1){40}}

\put(20,60){\line(1,-1){40}}

\put(40,80){\line(-1,-1){40}}
\put(40,80){\line(1,-1){40}}

\put(39,-4){$u$}
\put(39,82){$v$}

\put(40,0){\circle*{2}}
\put(20,20){\circle*{2}}
\put(60,20){\circle*{2}}
\put(0,40){\circle*{2}}
\put(20,40){\circle*{2}}
\put(40,40){\circle*{2}}
\put(60,40){\circle*{2}}
\put(80,40){\circle*{2}}
\put(20,60){\circle*{2}}
\put(60,60){\circle*{2}}
\put(40,80){\circle*{2}}
\end{picture}\caption{Lemma~\ref{L:lattice-to-chain} does not extend to lattices of greater length.}
\label{F:cel}
\end{figure}

The result of Lemma~\ref{L:lattice-to-chain} does not extend to length four or more. The lattice of Figure~\ref{F:cel} is not a congruence-preserving extension of any chain with extremities~$u$ and $v$.

\begin{lemma}\label{L:semnotliftable}
Let $n\ge 4$ be an integer, let~$\cV$ be a finitely generated variety of lattices such that $M_n\not\in\cV$. If $\lh(K)\le 3$ for each simple lattice~$K$ of~$\cV$, then $\crit{\cM_n^{0,1}}{\cV}\le\aleph_2$.
\end{lemma}

\begin{proof}
We consider the diagram $\vec A$ introduced just before Lemma~\ref{L:cong-for-Mn}. Assume that $\crit{\cM_n^{0,1}}{\cV}>\aleph_2$. As $M_n\in\cM_n^{0,1}$, $\vec A$ is a diagram of~$\cM_n^{0,1}$ indexed by~$I_n$. By Lemma~\ref{L:lifting-diagram-idexed-by-I}, the diagram $\Conc\circ\vec A$ has a congruence-lifting $\vec B=(B_P,g_{P,Q})_{P\subseteq Q\text{ in }I_n}$ in~$\cV$. As $\Con B_{\un}\cong 2$, the lattice $B_{\un}$ is simple, thus, by assumption on~$\cV$, $\lh(B_{\un})\le 3$, and so $\lh(B_{\set{x}})\le 3$, for all $x\in \un$. The lattice $B_\emptyset$ is simple, so, taking a sublattice, we can assume that $B_\emptyset=\set{u,v}$, with $u<v$. By Lemma~\ref{L:lattice-to-chain}, we can assume that $B_{\set{x}}$ is a chain of length two, for each $x\in\un$. So by Lemma~\ref{L:lifting-imply-Mn}, $M_n$ is a sublattice of $B_{\un}$, and so $M_n\in\cV$, a contradiction.
\end{proof}

\begin{theorem}
Let~$\cV$ be a finitely generated variety of modular lattices and~$\cW$ be finitely generated variety of lattices. Let $n\ge 3$ be an integer such that $M_n\in\cV - \cW$. If $\lh(K)\le 3$ for each simple $K\in\cV$, then $\crit{\cV}{\cW}\le\aleph_2$. Moreover if either $\lh(K)\le 2$ for each simple $K\in\cV$ and $M_3\in\cW$ or $\lh(K)\le 3$ for each simple $K\in\cV$ and $\Sub F^3\in\cW$ for some field~$F$, then $\crit{\cV}{\cW}=\aleph_2$.
\end{theorem}

\begin{proof}
By Lemma~\ref{L:semnotliftable}, $\crit{\cV}{\cW}\le\aleph_2$.

Assume that $\lh(K)\le 2$ for each simple $K\in\cV$ and $M_3\in\cW$. As $\Sub \FF_2^2\cong M_3\in\cW$, it follows from Theorem~\ref{T:crit-ge-aleph2} that $\crit{\cV}{\cW}\ge\aleph_2$.

Assume that $\lh(K)\le 3$ for each simple $K\in\cV$ and $\Sub F^3\in\cW$ for some field~$F$, it follows from Theorem~\ref{T:crit-ge-aleph2} that $\crit{\cV}{\cW}\ge\aleph_2$.
\end{proof}

Similarly we obtain the following critical points.

\begin{corollary}
The following equalities hold
\begin{align*}
\crit{\cM_n}{\cM_{m,m}}&=\aleph_2;\\
\crit{\cM_n^{0,1}}{\cM_{m,m}}&=\aleph_2;\\
\crit{\cM_n^{0,1}}{\cM_{m,m}^{0,1}}&=\aleph_2;\\
\crit{\cM_n}{\cM_{m,m}^{0}} &=\aleph_2;\\
\crit{\cM_n}{\cM_{m}^{0}} &=\aleph_2, &\text{for all $n,m$ with $3\le m<n\le\omega$.}
\end{align*}
\end{corollary}

\begin{proof}
Let $n'\le n$ be an integer such that $m<n'<\omega$. As $M_{n'}\not\in\cM_{m,m}$, it follows from Lemma~\ref{L:semnotliftable} that $\crit{\cM_{n'}^{0,1}}{\cM_{m,m}}\le\aleph_2$, thus:
\begin{equation}\label{Eq1:corcritMnMmm}
\crit{\cM_{n}^{0,1}}{\cM_{m,m}}\le\aleph_2.
\end{equation}
Moreover $M_3\in \cM_{m,m}$, simple lattices of~$\cM_{m,m}$ are of length at most $3$, and finitely generated lattices of~$\cM_n$ have finite length (and are even finite). Thus, by Theorem~\ref{T:crit-ge-aleph2}
\begin{equation}\label{Eq2:corcritMnMmm}
\crit{\cM_n}{\cM_{m,m}^0}\ge\aleph_2.
\end{equation}
Similarly:
\begin{equation}\label{Eq3:corcritMnMmm}
\crit{\cM_n^{0,1}}{\cM_{m,m}^{0,1}}\ge\aleph_2.
\end{equation}
All the desired equalities are immediate consequences of \eqref{Eq1:corcritMnMmm}, \eqref{Eq2:corcritMnMmm}, and \eqref{Eq3:corcritMnMmm}.
\end{proof}

As an immediate consequence we obtain:
\begin{corollary}
$\crit{\cM_{4,3}}{\cM_{3,3}}\le\aleph_2$.
\end{corollary}
This question was suggested by M. Plo\v s\v cica.

\begin{lemma}\label{L:crit-Mn-in-a-subF3}
Let~$F$ be field. Then $M_n\in\Var(\Sub F^3)$ if and only if $n\le 1 +\card F$.
\end{lemma}

\begin{proof}
If~$F$ is infinite then the result is obvious. So we can assume that~$F$ is finite.

If $n\le 1 +\card F$, then $M_n$ is a sublattice of $M_{1+\card F}\cong\Sub F^2\in \Var(\Sub F^3)$, thus $M_n\in \Var(\Sub F^3)$.

Now assume that $M_n\in \Var(\Sub F^3)$. By J\'onsson's Lemma, $M_n$ is a homomorphic image of a sublattice of $\Sub F^3$. As $M_n$ satisfies Whitman's condition, it follows from the Davey-Sands Theorem \cite[Theorem~1]{DaveySands} that $M_n$ is projective in the class of all finite lattices. Therefore, as $\Sub F^3$ is finite, $M_n$ is a sublattice of $\Sub F^3$. Thus there exist distinct subspaces $A,B,V_1,V_2,\dots,V_n$ of $F^3$ such that $V_i\cap V_j=A$ and $V_i + V_j=B$, for all $1\le i<j\le n$. Let $i,j,k$ distinct. Then:
\[
\dim V_i + \dim V_j = \dim B + \dim A = \dim V_i + \dim V_k.
\]
Thus $\dim V_j=\dim V_k$. But $\dim A<\dim V_1<\dim B\le \dim F^3=3$. If $\dim A=1$, then $M_n$ is isomorphic to $\set{A/A,V_1/A,\dots,V_n/A,B/A}$ which is a sublattice of $\Sub(B/A)$, with $\dim B/A=2$. If $\dim A=0$, then:
\[\dim B=\dim(V_1\oplus V_2)=\dim V_1 + \dim V_2= 2\cdot \dim V_1.\]
Thus $\dim B$ is even, moreover $\dim B\le 3$, hence $\dim B=2$.

In both cases $M_n$ is a sublattice of $\Sub E$ for some~$F$-vector space $E$ of dimension two. But $\Sub E\cong M_{1+\card F}$, thus $n\le 1+\card F$.
\end{proof}

\begin{corollary}
Let~$F$ be a finite field and let $n>1+\card F$. Then:
\begin{align*}
\crit{\cM_n}{\Var(\Sub F^3)}&=\aleph_2;\\
\crit{\cM_n}{\Var_0(\Sub F^3)}&=\aleph_2;\\
\crit{\cM_n^{0,1}}{\Var(\Sub F^3)}&=\aleph_2;\\
\crit{\cM_n^{0,1}}{\Var_{0,1}(\Sub F^3)}&=\aleph_2.
\end{align*}
\end{corollary}

\begin{proof}
By Lemma~\ref{L:crit-Mn-in-a-subF3}, $M_n\not\in \Var(\Sub F^3)$, moreover simple lattices of $\Var(\Sub F^3)$ are of length at most three. Thus, by Lemma~\ref{L:semnotliftable}:
\begin{equation}\label{Eq1:corcritMnfield}
\crit{\cM_n^{0,1}}{\Var(\Sub F^3)}\le\aleph_2.
\end{equation}

As each simple lattice of~$\cM_n$ is of length at most two, it follows from Theorem~\ref{T:crit-ge-aleph2} that
\begin{equation}\label{Eq2:corcritMnfield}
\crit{\cM_n}{\Var_0(\Sub F^n)}\ge\aleph_2, \quad\text{and}\crit{\cM_n^{0,1}}{\Var_{0,1}(\Sub F^n)}\ge\aleph_2.
\end{equation}
All the other desired equalities are consequences of \eqref{Eq1:corcritMnfield}, \eqref{Eq2:corcritMnfield}.
\end{proof}

\begin{corollary}
Let~$F$ and~$K$ be finite fields. If $\card F>\card K$ then:
\begin{align*}
\crit{\Var(\Sub F^3)}{\Var(\Sub K^3)}&=\aleph_2;\\
\crit{\Var(\Sub F^3)}{\Var_0(\Sub K^3)}&=\aleph_2;\\
\crit{\Var_{0,1}(\Sub F^3)}{\Var(\Sub K^3)}&=\aleph_2;\\
\crit{\Var_{0,1}(\Sub F^3)}{\Var_{0,1}(\Sub K^3)}&=\aleph_2.
\end{align*}
\end{corollary}

\begin{proof}
By Lemma~\ref{L:crit-Mn-in-a-subF3}, $M_{1+\card F}\not\in \Var(\Sub K^3)$, moreover simple lattices of $\Var(\Sub K^3)$ are of length at most three. Thus, by Lemma~\ref{L:semnotliftable}:
\begin{equation}\label{Eq1:corcritfieldfield}
\crit{\Var_{0,1}(\Sub F^3)}{\Var(\Sub K^3)}\le\aleph_2.
\end{equation}

As each simple lattice of $\Var(\Sub F^3)$ is of length at most three, it follows from Theorem~\ref{T:crit-ge-aleph2} that:
\begin{align}
\crit{\Var(\Sub F^3)}{\Var_0(\Sub K^n)}\ge\aleph_2,\label{Eq2:corcritfieldfield}\\
\crit{\Var_{0,1}(\Sub F^3)}{\Var_{0,1}(\Sub K^n)}\ge\aleph_2. \label{Eq3:corcritfieldfield}
\end{align}
All the other desired equalities are consequences of \eqref{Eq1:corcritfieldfield}, \eqref{Eq2:corcritfieldfield}, \eqref{Eq3:corcritfieldfield}.
\end{proof}

\begin{lemma}\label{L:critMomega}
Let~$\cV$ be a finitely generated variety of lattices \textup(resp., a finitely generated variety of lattices with~$0$\textup), let $m\ge 2$ an integer. Assume that for each simple lattice~$K$ of~$\cV$, there do not exist $b_0,b_1,\dots,b_{m-1}>u$ in~$K$ such that $b_i\wedge b_j=u$ \textup(resp., $b_0,b_1,\dots,b_{m-1}>0$ such that $b_i\wedge b_j=0$\textup), for all $0\le i<j\le m-1$. Then $\crit{\cM_{2m-1}^{0,1}}{\cV}\le\aleph_2$.
\end{lemma}

\begin{proof}
Set $n=2m-1\ge 3$. Let $\vec A=(A_P,f_{P,Q})_{P\subseteq Q\text{ in }I_n}$ be the direct system of~$\cM_n^{0,1}$ introduced just before Lemma~\ref{L:cong-for-Mn}. Assume that $\crit{\cM_{n}^{0,1}}{\cV}>\aleph_2$. By Lemma~\ref{L:lifting-diagram-idexed-by-I}, there exists a congruence-lifting $\vec B=(B_P,g_{P,Q})_{P\subseteq Q\text{ in }I_n}$ of $\Conc\circ\vec A$ in~$\cV$. Let $\vec \xi=(\xi_P)_{P\in I_n}\colon\Conc\circ\vec A\to\Conc\circ\vec B$ be a natural equivalence. Taking a sublattice of $B_\emptyset$, we can assume that $B_\emptyset$ is a chain $u<v$. Moreover, as the map $f_{P,Q}$ is an inclusion map, we can assume that $g_{P,Q}$ is an inclusion map, for all $P\subseteq Q$ in~$I_n$.

Let $x\in\un$. By Lemma~\ref{L:cong-for-Mn}, $\Theta_{B_{\set{x}}} (u,v)$ is the largest congruence of $B_{\set{x}}$. Thus:
\[
\Theta_{B_{\set{x}}} (u,v)=\xi_{\set{x}}(\Theta_{A_{\set{x}}} (0,a_x)) \vee \xi_{\set{x}}(\Theta_{A_{\set{x}}} (a_x,1))
\]
Therefore there exist $t_0^x=u<t_1^x<\dots<t_{r+1}^x=v$ in $B_{\set{x}}$ such that, for all $0\le i\le r$:
\[
\text{either }(t_{i}^x,t_{i+1}^x)\in \xi_{\set{x}}(\Theta_{A_{\set{x}}} (0,a_x))\text{ or }(t_{i}^x,t_{i+1}^x)\in \xi_{\set{x}}(\Theta_{A_{\set{x}}} (a_x,1))
\]
Set $b_x=t_1^x$. Put:
\[
X'=\setm{x\in \un}{\Theta_{B_{\set{x}}}(u,b_x) = \xi_{\set{x}}(\Theta_{A_{\set{x}}}(0,a_x))}\]
\[X''=\setm{x\in \un}{\Theta_{B_{\set{x}}}(u,b_x) = \xi_{\set{x}}(\Theta_{A_{\set{x}}}(a_x,1))}.\]
By symmetry we can assume that $\card X'\ge\card X''$ (we can replace the diagram $\vec A$ by its dual if required). As $\un=X'\cup X''$ and $\card\un=n=2m-1$, $\card X'\ge m$. Let $x,y\in X'$ distinct, it follows from Lemma~\ref{L:cong-for-Mn}(2) that $b_x\wedge b_y=u$. So we obtain a family of elements $(b_x)_{x\in X'}$ greater than $u$ such that $b_x\wedge b_y=u$ (resp., $b_x\wedge b_y=u=0$) for all $x\not=y$ in $X'$, a contradiction.
\end{proof}

With a similar proof using Lemma~\ref{L:lifting-diagram-idexed-by-I-nonFG} instead of Lemma~\ref{L:lifting-diagram-idexed-by-I} we obtain the following lemma.
\begin{lemma}\label{L:critMomega2}
Let~$\cV$ be a variety of lattices \textup(resp., a variety of lattices with~$0$\textup), let $m\ge 2$ an integer. Assume that for each simple lattice~$K$ of~$\cV$, there do not exist $b_0,b_1,\dots,b_{m-1}>u$ in~$K$ such that $b_i\wedge b_j=u$ \textup(resp., $b_0,b_1,\dots,b_{m-1}>0$ such that $b_i\wedge b_j=0$\textup), for all $0\le i<j\le m-1$. Then $\crit{\cM_{2m-1}^{0,1}}{\cV}\le\aleph_3$.
\end{lemma}

\begin{theorem}
Let~$\cV$ be either a finitely generated variety of lattices or a finitely generated variety of lattices with~$0$. If $M_3\in\cV$ then:
\begin{align*}
\crit{\cM_\omega}{\cV} &=\aleph_2;\\
\crit{\cM_\omega^0}{\cV}&=\aleph_2.
\end{align*}
Let~$\cV$ be a finitely generated variety of bounded lattices. If $M_3\in\cV$ then:
\begin{align*}
\crit{\cM_\omega^{0,1}}{\cV}=\aleph_2.
\end{align*}
\end{theorem}

\begin{proof}
Let~$\cV$ be a finitely generated variety of lattices, let~$m$ be the maximal cardinality of a simple lattice of~$\cV$. Thus the assumptions of Lemma~\ref{L:critMomega} are satisfied, so \emph{a fortiori} $\crit{\cM_{2m-1}^{0,1}}{\cV}\le\aleph_2$, and so $\crit{\cM_{\omega}^{0,1}}{\cV}\le\aleph_2$.

Denote by $\FF_2$ the two-element field. Let~$\cV$ be a variety of lattices with~$0$ (resp., with~$0$ and~$1$), such that $M_3\in\cV$. The variety~$\cM_{\omega}$ is locally finite, thus all finitely generated lattices of~$\cM_{\omega}$ are of finite length. Moreover all simple lattices of~$\cM_\omega$ have length at most two. Thus, by Theorem~\ref{T:crit-ge-aleph2}:
\[
\crit{\cM_{\omega}}{\Var_0(\Sub \FF_2^2)}\ge\aleph_2\text{ (resp., }\crit{\cM_{\omega}^{0,1}}{\Var_{0,1}(\Sub \FF_2^2})\ge\aleph_2).
\]
Moreover $\Sub \FF_2^2\cong M_3$, so $\crit{\cM_{\omega}}{\cV}\ge\aleph_2.$
\end{proof}

\section{Acknowledgment}
I thank Friedrich Wehrung whose help, advise, and careful reading of the paper, led to several improvements both in form and substance. I also thank Miroslav Plo\v{s}\v{c}ica for the stimulating discussions that led to the diagram~$\vec{A}$ introduced in Section~\ref{S:majoration-critpoint}. Partially supported by the institutional grant MSM 0021620839.

\end{document}